\documentclass[10pt]{amsart}
\usepackage{import}
\usepackage{./preamble/Pakete}
\usepackage{./preamble/Befehle}

\begin{document}

%!TEX root = ../main.tex

\title{An Intersection Matrix for Affine Hyperplane Arrangements}

\author{Jens Niklas Eberhardt}
\address{Mathematical Institute, University of Bonn
}
\email{mail@jenseberhardt.com}  % 'optional'
\author{Carl Mautner}
\address{
Department of Mathematics,
University of California Riverside
%900 University Ave.\\
%Riverside, CA 92521
}
\email{mautner@math.ucr.edu}
%\subjclass{?)}
%\keywords{?}
%\begin{abstract}
%?
%\end{abstract}
\maketitle

\begin{abstract}
For a real affine hyperplane arrangement, we define  an integer intersection matrix with a natural $q$-deformation related to the intersections of bounded chambers of the arrangement.
By connecting the integer matrix to a bilinear form of Schechtman--Varchenko, we show that there is a closed formula for its determinant that only depends on the combinatorics of the underlying matroid. We conjecture an analogous formula for its $q$-deformation. Our work also applies more generally in the setting of affine oriented matroids.

Additionally, we give a representation-theoretic interpretation of our $q$-intersection matrix using Braden--Licata--Proudfoot--Websters's hypertoric category $\mathcal{O}$ (or more generally Kowalenko--Mautner's category $\mathcal{O}$ for oriented matroid programs).  This paper is part of a broader program to categorify matroidal Schur algebras defined by Braden--Mautner.
\end{abstract}

%!TEX root = ../main.tex
\section{Introduction}

\subsection{An intersection matrix associated to hyperplane arrangements} Let $\cC$ be an arrangement of $n$ affine hyperplanes in $\R^r$. We assume that the hyperplanes are generically translated, so that any $d$ hyperplanes with nonempty intersection intersect in codimension $d$. We further assume that the associated central arrangement $\cC_0$ is essential.

We now explain how to associate an integer intersection matrix to $\cP_0(\cC)$, the set of bounded regions of the arrangement $\cC$ (meaning bounded connected components of the complement of the hyperplanes $\R^r - \bigcup_{H \in \cC} H$).

Note that our assumptions imply that the closure $\overline{A}$ of any bounded region $A$ is a simple polytope.  For a simple polytope $P$ of dimension $d$, let $f_0(P)$ be the number of vertices of $P$, and more generally for $i=0,1,\ldots,d$, let $f_i(P)$ denote the number of $i$-faces of $P$.  For regions $A,B$, let $d(A,B)$ be the number of hyperplanes separating them.

\begin{definition}
For bounded regions $A,B \in P_0(\cC)$, let $$S(A,B) = (-1)^{d(A,B)} \cdot f_0(\overline{A}\cap \overline{B}).$$
Fix an ordering on the set $\cP_0(\cC)$.  The \textit{intersection matrix} $S(\cC)$ is the matrix with rows and columns labelled by $\cP_0(\cC)$ and $(A,B)$-entry given by $S(A,B)$.
\end{definition}

\begin{remark}\label{rmk-HS}
For unimodular hyperplane arrangements, the matrix $S(\cC)$ appeared in~\cite{HS} where it was identified with an intersection form in the cohomology of the corresponding smooth hypertoric variety, see Section \ref{sec-interpret}.
\end{remark}

\begin{example}\label{ex-standard} Consider the following two arrangements $\mathcal{C}$ and $\mathcal{C'}$ with two bounded chambers that differ by a translation of one of the hyperplanes.
%\begin{figure}[ht]
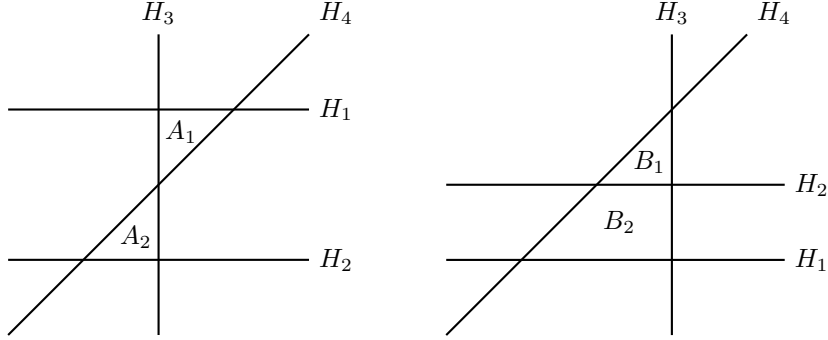
\[
    \centering
    \begin{minipage}{0.48\textwidth}
    \centering
    \begin{tikzpicture}
        % Define lines
        \draw[thick] (-2,1) -- (2,1) node[right] {$H_1$}; % Horizontal line 1
        \draw[thick] (-2,-1) -- (2,-1) node[right] {$H_2$}; % Horizontal line 2
        \draw[thick] (0,-2) -- (0,2) node[above] {$H_3$}; % Vertical line
        \draw[thick] (-2,-2) -- (2,2) node[above right] {$H_4$}; % 135-degree line
    
        % Label chambers
        \node at (0.3, 0.7) {$A_1$};
        \node at (-0.3, -0.7) {$A_2$};
    \end{tikzpicture}
    \end{minipage}\hfill
    \begin{minipage}{0.48\textwidth}
    \centering
    \begin{tikzpicture}
        % Define lines
        \draw[thick] (-3,-2) -- (1.5,-2) node[right] {$H_1$}; % Lower horizontal line
        \draw[thick] (-3,-1) -- (1.5,-1) node[right] {$H_2$}; % Upper horizontal line
        \draw[thick] (0,-3) -- (0,1) node[above] {$H_3$}; % Vertical line
        \draw[thick] (-3,-3) -- (1,1) node[above right] {$H_4$}; % Shifted 135-degree line
    
        % Label chambers
        \node at (-.3, -.7) {$B_1$};
        \node at (-0.7, -1.5) {$B_2$};
    \end{tikzpicture}
    \end{minipage}
%    \caption{Two arrangements $\mathcal{C}$ (left) and $\mathcal{C'}$ (right) with the same underlying matroid.}
%    \label{fig:thetwostandardarrangements}
\]%    \end{figure}
The corresponding intersection matrices are
\[
S(\mathcal{C}) = \begin{pmatrix}
3 & 1 \\
1 & 3
\end{pmatrix}, \quad
S(\mathcal{C}') = \begin{pmatrix}
3 & -2 \\
-2 & 4
\end{pmatrix}.
\]

\end{example}

A first objective of this paper is to relate the matrix $S$ to a bilinear form defined by Schechtman--Varchenko~\cite{SV} for hyperplane arrangements.  Then using a formula of Schechtman--Varchenko it follows that the determinant of $S$ can be described in terms of the combinatorics of the underlying matroid $M$ associated to $\cC$ defined by normal vectors to the affine hyperplanes.

To state the determinant formula we need the following notation.  Let $I$ be the ground set of the matroid $M$.  Denote by $\cF$ the set of coloop-free flats and by $\mu$ the M\"obius function of the matroid $M$.  For any flat $K$, let $\mu^+(K) = (-1)^{r(K)}\mu(\emptyset,K)$ be the unsigned M\"obius function, which is always positive and equal to the number of nbc bases of $K$.  Crapo's beta invariant $\beta(K)$ for a flat $K$ is defined as
\[ \beta(K) := (-1)^{r(K)} \sum_{K': K'\leq K} \mu(K') r(K'),\]
where $K'$ runs over all flats contained in $K$.  The beta invariant can be computed by using the following formula, which holds for any $e\in I$ that is neither a loop nor coloop,
\[ \beta(M) = \beta(M-e)+\beta(M/e),\]
together with the facts that $\beta(M)=0$ if $M$ is a loop or is not connected and $\beta(M)=1$ if $M$ is a coloop.  In particular it follows that, for any matroid $M$, $\beta(M)\geq 0$.
\begin{theorem}[Determinant formula]\label{thm-det}
The determinant of the matrix $S(\cC)$ is 
\[ \det S(\cC) = \prod_{K \in \cF \setminus\{I\}} |I\setminus K|^{\beta(M/K) \mu^+(K)} .\]
\end{theorem}

\begin{example} In Example~\ref{ex-standard}, the proper coloop-free flats in the underlying matroid $M$ are $\emptyset$ and  $F=\{H_1,H_2\}.$ Indeed, we find that
$$\det(\cC) = \det(\cC') = 8 = 4^1 \cdot 2^1 = |I|^{\beta(M)\mu^+(\emptyset)}|I\setminus F|^{\beta(M/F)\mu^+(F)}.$$
\end{example}

\subsection{A $q$-intersection matrix} The second objective of this paper is to introduce a $q$-version of the intersection matrix and a conjectural determinant formula.  

For a simple polytope $P$ of dimension $d$, the $h$-polynomial of $P$ is defined as $h(x) = f(x-1)$, where $f(x) = f_0+f_1 x + f_2 x^2 +\ldots + f_d x^d$.  In particular, note that $h(1) = f(0) = f_0$. 

\begin{definition}
For bounded regions $A,B \in P_0(\cC)$, let $$S_q(A,B) = (-q)^{d(A,B)} \cdot h(\overline{A}\cap \overline{B},q^2),$$
where $h(P,q^2)$ is the $h$-polynomial of $P$ evaluated at $q^2$.

Fix an ordering on the set $\cP_0(\cC)$.  The $q$-\textit{intersection matrix} $S_q(\cC)$ is the matrix with rows and columns labelled by $\cP_0(\cC)$ and $(A,B)$-entry given by $S_q(A,B)$.
\end{definition}

Note that if we set $q=1$ in the expression for $S_q(A,B)$, we recover $S(A,B) = (-1)^{d(A,B)} \cdot f_0(\overline{A}\cap \overline{B}),$ so $S_1(\cC) = S(\cC)$.  

In Theorem~\ref{EulerForm} we give an algebraic interpretation of the matrix $S_q(\cC)$.  Namely, we show that the $q$-matrix $S_q(\cC)$ appears as a Gram matrix for the Euler form on the Grothendieck group of hypertoric category $\cO$ defined by Braden--Licata--Proudfoot--Webster~\cite{GDKD,BLPWtorico}.

For a positive integer $n$, let $[n]_{q^2} = 1 + q^2 + \ldots + q^{2n-2}$.

\begin{Vermutung}\label{conj:mainconjecture}
The determinant of the matrix $S_q(\cC)$ is 
\[ \det S_q(\cC) = \prod_{K \in \cF \setminus\{I\}} [| I\setminus K |]_{q^2}^{\beta(M/K) \mu^+(K)} .\]
\end{Vermutung}

\begin{example} Returning to Example~\ref{ex-standard}, the $q$-intersection matrices are
\[
S_q(\mathcal{C}) = \begin{pmatrix}
1+q^2+q^4 & (-q)^2 \\
(-q)^2 & 1+q^2+q^4
\end{pmatrix}, \quad
S_q(\mathcal{C}') = \begin{pmatrix}
1+q^2+q^4 & -q(1+q^2) \\
-q(1+q^2) & 1+2q^2+q^4
\end{pmatrix}.
\]
It follows that
$$\det S_q(\mathcal{C})=\det S_q(\mathcal{C}')=(1+q^2+q^4+q^6)(q^2+1)=[4]_{q^2}[2]_{q^2}.$$
The coloop-free flats $\neq I=\{H_1,\dots,H_4\}$ in the underlying matroid $M$ are $\emptyset$ and  $F=\{H_1,H_2\}$ and so
$$ \prod_{K \in \cF \setminus\{I\}} [| I\setminus K |]_{q^2}^{\beta(M/K) \mu^+(K)} = [|I|]^{\beta(M)\mu^+(\emptyset)}_{q^2}[|I\setminus F|]^{\beta(M/F)\mu^+(F)}_{q^2}=[4]_{q^2}[2]_{q^2}$$
in agreement with our conjecture.
\end{example}

\subsection{Motivation from matroidal Schur algebras}
In~\cite{BMmatroid}, Braden and the second author introduced matroidal Schur algebras, a class of cellular (in fact quasi-hereditary) algebras associated to any matroid $M$, and showed that the canonical bilinear forms associated to their cell modules can be identified with Schechtman--Varchenko's (or more generally Brylawski--Varchenko's~\cite{BV}) bilinear forms on flag space for the dual matroid $M^*$.  From the proof of Theorem~\ref{thm-det} it follows that for any hyperplane arrangement $\cC$ as above, the matrix $S(\cC)$ is a Gram matrix for the bilinear form on a cell module for the matroidal Schur algebra associated to the matroid underlying $\cC_0$.

We expect that for the matroid represented by $\cC_0$ there is a natural $q$-analogue of the matroidal Schur algebras for which $S_q(\cC)$ is a Gram matrix for the bilinear form on a cell module and that this $q$-matroidal Schur algebra is categorified by hypertoric Harish-Chandra bimodules.

For general non-representable matroids, it is not clear to us how to define an analogue of $S_q(\cC)$.  However, in the next section we explain how our definitions and results do naturally extend to the setting of affine oriented matroids.  This gives hope that one could define $q$-matroidal Schur algebras at least for orientable matroids.

\subsection{Intersection matrices for affine oriented matroids}  Motivated by the expectation that there is a $q$-analogue of the matroidal Schur algebras of~\cite{BMmatroid}, in this section we extend our theorem and conjecture to the broader setting of affine oriented matroids.  A reader who is only interested in hyperplane arrangements may skip this section.  We will assume knowledge of the basics of oriented matroid theory and refer to the book~\cite{OMbook} for definitions.

In this context, the role of the central arrangement $\cC_0$ is played by an oriented matroid $\cM$, while the role of the affine arrangement $\cC$ is played by an affine oriented matroid $(\tilde\cM,g)$.  We assume that $\tilde\cM/g = \cM$ and that $(\tilde\cM,g)$ is generic in the sense that for any cocircuit $Y$ of $\tilde\cM$, if $|z(Y)|>d$, then $Y(g)=0$.  Note that every central arrangement $\cC_0$ and any arrangement $\cC$ obtained from generic translation of the hyperplanes as above naturally gives rise to an underlying oriented matroid $\cM$ and affine oriented matroid $(\tilde\cM,g)$, generic in the sense described (see ~\cite[Example 2.4]{KM}).

Let $\cL$ denote the set of covectors of $\tilde\cM$ and let $\cL^+ = \{X \in \cL ~|~X(g)=+\}$ be the affine face lattice of $\tilde\cM$.  Let $\cL^{++}$ be the bounded complex, defined as the set of covectors, all of the faces of which are in $\cL^+$.  The maximal elements of $\cL^+$ are the topes of $\tilde\cM$.  Let $\cP_0(\tilde\cM)$ denote the bounded topes of $\tilde\cM$, or equivalently the maximal elements of $\cL^{++}$.  By our assumption on $\tilde\cM$, for any $A \in \cP_0(\tilde\cM)$, the opposite poset of all faces of $A$ in $\cL$, known as the Las Vergnas face lattice $\mathfrak{F}_{lv}(A)$ of $A$, is the face lattice of a pure simplicial complex, whose vertices correspond to the facets of $Y$ and whose maximal simplices correspond to the cocircuit faces of $Y$.

Like we did for polytopes, for any bounded covector $Y$ of rank $d+1$, let $f_0(Y)$ be the number of cocircuit faces of $Y$, and more generally for $i=0,1,\ldots,d$, let $f_i(Y)$ denote the number of rank $i+1$-faces of $Y$ and $h(Y,x) = f(Y,x-1)$, where $f(Y,x) = f_0(Y)+f_1(Y) x + f_2(Y) x^2 +\ldots + f_d(Y) x^d$.  For topes $A,B$, let $d(A,B)$ be the number of $i \in I$ such that $A(i) \neq B(i)$.

\begin{definition}
For topes $A,B \in P_0(\tilde\cM)$, let $$S(A,B) = (-1)^{d(A,B)} \cdot f_0(A\wedge B)$$ and $$S_q(A,B) = (-q)^{d(A,B)} \cdot h(A\wedge B,q^2),$$
Fix an ordering on the set $\cP_0(\tilde\cM)$ and let $S(\tilde\cM)$ (resp. $S_q(\tilde\cM)$) be the matrix with rows and columns labelled by $\cP_0(\tilde\cM)$ and $(A,B)$-entry given by $S(A,B)$ (resp. $S_q(A,B)$).
\end{definition}

\begin{remark}
In the setting of affine oriented matroids, Theorem~\ref{EulerForm} also gives an algebraic interpretation of this matrix.  Namely, $S_q(\tilde\cM)$ appears as a Gram matrix for the Euler form on the Grothendieck group of category $\cO$ for an oriented matroid program defined in~\cite{KM}, which generalizes hypertoric category $\cO$.
\end{remark}

Let $M$ be the underlying (unoriented) matroid of $\cM$ and again let $\cF$ denote the set of coloop-free flats of $M$.  Then we have the following extension of Theorem~\ref{thm-det}:

\begin{theorem}\label{thm-det-om}
The determinant of the matrix $S(\tilde\cM)$ is 
\[ \det S(\tilde\cM) = \prod_{K \in \cF \setminus\{I\}} |I\setminus K|^{\beta(M/K) \mu^+(K)} .\]
\end{theorem}

We expect that Conjecture~\ref{conj:mainconjecture} also holds in this more general setting:

\begin{Vermutung}\label{conj:om}
The determinant of the matrix $S_q(\cC)$ is 
\[ \det S_q(\tilde\cM) = \prod_{K \in \cF \setminus\{I\}} [| I\setminus K |]_{q^2}^{\beta(M/K) \mu^+(K)} .\]
\end{Vermutung}

\subsection{Relation to other work}  The bilinear forms of Schechtman--Varchenko~\cite{SV}, which play a role in this paper, were originally introduced in the context of complex hyperplane arrangements, where each hyperplane $H$ is equipped with a complex number `exponent' $a(H)$.  In particular, one could also obtain a $q$-analogue of the matrix $S(\cC)$ by setting $a(H)=q$ for all $H\in \cC$, but this does not recover the $q$-analogue we describe above.  

The Schechtman--Varchenko determinant formula is similar in flavor to (and in fact a `quasiclassical' version of) another bilinear form introduced by Varchenko~\cite{V} associated to a weighted real affine hyperplane arrangement.  Varchenko's form is defined on the vector space with basis parametrized by the set of all (not just bounded) regions of the arrangement and the pairing for two regions $A$ and $B$ is defined to be the product of the weights all hyperplanes that separate the regions.  It would be interesting to understand the relation between Varchenko's matrix and the intersection matrix considered here.  Note that we identify our intersection matrix with the Schechtman--Varchenko form for the Gale dual arrangement, so we expect a duality involved in any connection to the Varchenko form.  A more recent survey of Varchenko's work and proof of his determinant formula appears in Aguiar--Mahajan's book~\cite[Chapter 8]{AgM}.

In our proof of Theorem~\ref{thm-det} we describe a basis, parametrized by compact regions of $\cC$, of the flag space for the Gale dual arrangement $\cC_0^\perp$.  Bj\"orner--Wachs~\cite{BW} also define a basis of the same flag space, but parameterized by the set of chambers of $\cC_0^\perp$ bounded with respect to a generic covector.  Under Gale duality, the choice of a generic translation $\cC$ of $\cC_0$ corresponds to a choice of generic covector $\xi$ for $\cC_0^\perp$ and the set of compact regions of $\cC$ are in canonical bijection with the set of $\xi$-bounded regions in $\cC_0^\perp$.  It would also be interesting to compare these two bases.

As mentioned in Remark~\ref{rmk-HS}, the integer matrix $S$ for unimodular hyperplane arrangements appears in~\cite{HS}, where it arises from a cohomological intersection pairing.  The authors of that paper also show that the form is identified with a pairing coming from the independence complex of the underlying matroid and we note that their proof contains versions of Lemmas~\ref{lem-pairing},~\ref{lem-U}, and~\ref{lem-basis} below.

\subsection{Outline} Section~\ref{sec-proof} contains a proof of the main result.  In Section~\ref{sec-interpret} we describe how the matrices under consideration arise in the geometry of hypertoric varieties and show that they can be realized in representation theory as a Euler form on the Grothendieck group of hypertoric category $\cO$.  In Appendix~\ref{sec-examples} we conclude with examples illustrating the conjecture, particularly in the affine oriented matroid setting.

\subsection{Acknowledgements} The second author thanks Ethan Kowalenko and Maksim Maydanskiy, who were helpful in checking the V\'amos example.  The first author was supported by Deutsche Forschungsgemeinschaft (DFG), project number 45744154, Equivariant $K$-motives and Koszul duality. The second author was supported in part by NSF grant DMS-1802299.

%!TEX root = ../main.tex

\section{Proof of Theorem~\ref{thm-det} and~\ref{thm-det-om}}\label{sec-proof}

While Theorem~\ref{thm-det} can be viewed as a special case of Theorem~\ref{thm-det-om} (and so it would be enough to prove~\ref{thm-det-om}), we prove the two statements in parallel in the hopes that the paper is still accessible to those who are not familiar with the oriented matroid setting.

Let $\Lambda(I)$ be the exterior algebra over $\Z$ on the set $I$. Let $\cB = \Lambda^r_r \subset \Lambda(I)$ be the sublattice generated by monomials corresponding to bases of $M$.

In the hyperplane arrangement setting, fix an orientation of our underlying vector space $\R^r$ and a normal direction for each hyperplane.  In particular, each region $A$ of the arrangement is then described by a sign vector $A:I \to \{+,-\}$ encoding whether the region lies on the positive or negative side of each hyperplane.  For each basis $b \subset I$, let $e_b = i_1 \wedge \ldots \wedge i_r$, where $i_1,\ldots,i_r$ is an ordering of $b$ such that the orientation defined by the ordered list of corresponding normal vectors agrees with our fixed choice of orientation of $\R^r$.

In the affine oriented matroid setting, we fix a choice of basis orientation $\chi$ of $\cM$ (see~\cite[Def. 3.5.1]{OMbook}). For a basis $b$ of $M$, let $e_b = \chi(i_1,\ldots,i_r)\cdot i_1 \wedge \ldots \wedge i_r$, where $i_1, \ldots, i_r$ is an arbitrary ordering of $b$.  By the definition of $\chi$, $e_b$ is well-defined.

Note for hyperplane arrangements that there is a natural bijection between the set of bases of $M$ and the set of vertices of $\cC$, where a basis $b$ of $M$ corresponds to the intersection $Y_b = \bigcap_{i\in b} H_i$.   In the affine oriented matroid setting, this translates to a bijection between the set of bases of $M$ and the set of feasible cocircuits~\cite[Lemma 2.6]{KM} and we let $Y_b$ be the feasible cocircuit corresponding to a basis $b$.

We equip $\cB$ with the perfect pairing defined by $\langle e_b,e_{b'} \rangle = \delta_{b,b'}$.

Let $\cU \subset \cB = \Lambda^r_r$ be the kernel of the restriction of the standard differential $\partial: \Lambda(I) \to \Lambda(I)$ to $\Lambda^r_r$.  By~\cite[Theorem 4.1]{BMmatroid} we can identify $\cU$ and the flag space of Brylawski--Schechtman--Varchenko $\cF^{n-r}(M^*)$ for the dual matroid and under this identification the restriction of the perfect pairing $\langle~,~\rangle$ to $\cU$ agrees with the Brylawski--Schechtman--Varchenko bilinear form $S$~\cite{SV,BV} for $M^*$.

For each bounded region (or tope) $A \in \cP_0$, let
\[ \phi(A) = \sum_b \left(\prod_{i\in b} A(i)\right) \cdot e_b \in \cB,\]
where the sum runs over all bases $b$ such that $Y_b$ is a vertex (cocircuit face) of $A$ and $A(i)$ is the sign vector for $A$ evaluated at $i$.

\begin{lemma}\label{lem-pairing}
For any $A,B \in \cP_0$, there is an equality
\[ \langle \phi(A),\phi(B)\rangle = S(A,B).\]
\end{lemma}

\begin{proof}
The pairing $\langle \phi(A),\phi(B)\rangle = \langle \sum (\prod_{i\in b} A(i)) \cdot e_b, \sum (\prod_{i\in b'} B(i)) \cdot e_{b'}\rangle$ is equal to the sum over common vertices (cocircuit faces) $Y_b$ of $A$ and $B$, of $$\prod_{i\in b} A(i)\prod_{i\in b} B(i) = \prod_{i\in b} A(i)B(i).$$  Now if $b$ is a common vertex, then the sign vectors for $A$ and $B$ agree for all hyperplanes $H_j$ where $j \not\in b$ (as $A$ and $B$ can only be separated by hyperplanes in $b$).  Thus for any common vertex $b$, the sign $\prod_{i\in b} A(i)B(i) = \prod_{i\in I} A(i)B(i)= (-1)^{d(A,B)}$ only depends on the the number $d(A,B)$.  The number of common vertices is equal to $f_0(\overline{A}\cap \overline{B})$ and so we conclude that $$\langle \phi(A),\phi(B)\rangle =  (-1)^{d(A,B)} \cdot f_0(\overline{A}\cap \overline{B}) = S(A,B),$$
as we wished to show. \end{proof}

\begin{lemma}\label{lem-U}
For each $A \in \cP_0$, $\phi(A) \in \cU$.
\end{lemma}

\begin{proof}
We need to show that $\partial (\phi(A)) = 0$.  By definition, $\phi(A) = \sum_b (\prod_{i\in b} A(i)) \cdot e_b$ and so $\partial (\phi(A)) = \sum_b (\prod_{i\in b} A(i)) \cdot  \partial e_b$.  Now a monomial $e_J$ can only appear in $\partial e_b$ with non-zero coefficient if $J=b\setminus {i}$ for some $i$.  The intersection of the hyperplanes corresponding to elements of $J$ is a line (or rank 2 pseudosphere) and so intersects the closure of the polytope $A$ in an edge.  Thus there are only two bases corresponding to vertices of $A$ such that $e_J$ can appear in the differential, say $b_1 = J \sqcup \{i_1\}$ and $b_2 = J \sqcup \{i_2\}$.

Fix an ordering of $J = \{j_1,\ldots,j_{r-1}\}$ such that the orientation corresponding to the ordered list $i_1,j_1,\ldots,j_{r-1}$ agrees with the chosen orientation of $\R^r$. Then $e_{b_1}= i_1 \wedge j_1 \wedge \ldots \wedge j_{r-1}$. Note that $A(i_1) =\pm A(i_2)$ if and only if $e_{b_2} = \mp i_2 \wedge j_1 \wedge \ldots \wedge j_{r-1}$.  In other words, $e_{b_2}= -A(i_1)A(i_2) i_2 \wedge j_1 \wedge \ldots \wedge j_{r-1}$.

We conclude that the only terms of $\phi(A) = (\prod_{i\in {b}} A(i)) \cdot  \partial e_b$ whose boundary could contribute a non-zero multiple of $j_1 \wedge \ldots \wedge j_{r-1}$ are $(\prod_{i\in {b_1}} A(i)) \cdot  \partial e_{b_1}$ and $(\prod_{i\in {b_2}} A(i)) \cdot  \partial e_{b_2}$. By our choice of ordering, the multiplicity of $j_1 \wedge \ldots \wedge j_{r-1}$ in $(\prod_{i\in {b_1}} A(i)) \cdot \partial e_{b_1}$ is given by $A(i_1) (\prod_{j \in J} A(j))$ and its multiplicity in $(\prod_{i\in {b_2}} A(i)) \cdot \partial e_{b_2}$ is given by $$A(i_2) (\prod_{j \in J} A(j)) (-A(i_1)A(i_2))= -A(i_1) (\prod_{j \in J} A(j)).$$  Thus there is cancellation and  $\partial (\phi(A))=0$.

In the oriented matroid setting the same cancellation follows from the dual pivoting property (PV*) for basis orientations (see~\cite[p. 125]{OMbook}). \end{proof}

\begin{lemma}\label{lem-basis}
The set $\{\phi(A)~|~A \in \cP_0\}$ forms a $\Z$-basis of $\cU$.
\end{lemma}

\begin{proof}
We recall that $\cU = \Ker(\partial) = \tilde{H}_{r-1}(IN(M))$ is the top reduced homology of the independence complex of $M$, which by~\cite[Theorem 7.8.1]{Bj} is isomorphic to $\Z^{\mu^+(M^*)}$.  Note that as $\cU$ is the kernel of a map between free modules, $\cU$ is a direct summand of the lattice $\cB$.

We will first show that the cardinality of $\cP_0$ is equal to ${\mu^+(M^*)}$, the rank of $\cU$. For a hyperplane arrangement, by Gale duality, the set $\cP_0$ of compact regions in $\cC$ is in natural bijection with the regions in $\cC_0^\perp$ that are bounded for the corresponding generic covector.  By~\cite[Theorem 3.2]{GZ}, the number of such regions is equal to $\mu^+(M^*)$.

More generally, duality for oriented matroid programs~\cite[Corollary 10.1.11]{OMbook} gives a natural bijection between $\cP_0$ and the set of topes of $\tilde\cM^\vee$ bounded with respect to $g$.  By our genericity condition on $g$, the (unoriented) matroid underlying $\tilde\cM^\vee$ is the free single-element extension of $M^*$ (in the sense of~\cite[Section 7.2 and Exercise 7.2.3]{Ox}).  It follows that the flats of the matroid underlying $\tilde\cM^\vee$ are just the flats of $M^\vee$, but with $g$ added to the flat $I$.  A version of Zaslavsky's theorem on the number of bounded regions~\cite[Proof of Theorem 4.6.5]{OMbook} then gives:
\[ \#\cP_0 =  (-1)^{r-1} \sum_{B} \mu(\emptyset,B), \]
where the sum runs over all flats of the matroid underlying $\tilde\cM^\vee$ that do not contain $g$ or in other words all proper flats of $M^*$.  We conclude that $\cP_0 = \mu^+(M^*)$.

By Lemma~\ref{lem-U} and the fact that $\cU$ is a direct summand of $\cB$ of rank equal to the number of bounded regions (or bounded topes), it suffices to show that $\{\phi(A)~|~A \in \cP_0\}$ is part of a $\Z$-basis for $\cB$.

Let $\xi \in (\R^r)^*$ be a generic choice of covector (generic in the sense that $\xi$ is nonconstant when restricted to any positive dimensional face of $\cC$).  Analogously in the oriented matroid setting, we choose a generic Euclidean single-element extension $(\tilde\cM,g,f)$ of the affine oriented matroid.  The choice of $\xi$ (or $f$) allows us to define the set $\cP$ of $\xi$-bounded regions (or $f$-bounded topes).   Optimizing the function $\xi$ on each $\xi$-bounded region gives a natural bijection $\mu$ between the set $\Bas$ of bases of $M$ and the set $\cP$.  More generally, in the oriented matroid setting there is a natural bijection~\cite[Cor. 2.15]{KM} $\mu$ between $\Bas$ and the set $\cP$, which takes a basis $b$ to the unique tope whose optimal cocircuit $Y_b$ has zero set equal to $b$.

We extend the set $\{\phi(A)~|~A \in \cP_0\}$ by extending the map $\phi$ to all of $\cP$.  Namely, for $A \in \cP$, we let
\[ \phi(A) = \sum_b \left(\prod_{i\in b} A(i)\right) \cdot e_b \in \cB,\]
where the sum runs over all bases $b$ such that $Y_b$ is a vertex (feasible cocircuit face) of $A$ and $A(i)$ is the sign vector for $A$ evaluated at $i$.

The same argument as in Lemma~\ref{lem-pairing} shows that more generally for any $A,B \in \cP$,
\[ \langle \phi(A),\phi(B)\rangle = (-1)^{d(A,B)} \cdot \#\{\text{{feasible vertices of }} A \wedge B \}.\]

Let $\mathbf{e}_b = (\prod_{i\in b} (\mu(b))(i)) e_b.$  Note that $\{\mathbf{e}_b~|~b\in \Bas\}$ is a $\Z$-basis for $\cB$ and that
\[ \phi(A) = \sum_b \left(\prod_{i\in b} A(i)\right) \cdot e_b = \sum_b \left(\prod_{i\in b} A(i)\cdot (\mu(b))(i)\right) \cdot \mathbf{e}_b \]
\[  = \sum_b (-1)^{d(A,\mu(b))} \cdot \mathbf{e}_b \]
because $A(i) = \mu(b)(i)$ for all $i \not\in b$.

As the set $\{\mathbf{e}_b~|~b \in \Bas\}$ is a $\Z$-basis for $\cB$, to show that $\{\phi(A)~|~A \in \cP\}$ is a basis for $\cB$, it suffices to show that the matrix $y$ expressing $\phi(A)$ for each $A \in \cP$ in terms of the $\mathbf{e}_b$'s is invertible.

Note that $y$ has entries:
\[y_{A,B} = \begin{cases}
(-1)^{d(A,B)} & \mathrm{if~}Y_{\mu^{-1}(\b)} \mathrm{~is~a~face~of~}A \\
0 & \mathrm{otherwise},
\end{cases}\]

and $y$ has a natural $q$-analogue $Y$ with entries:
\[Y_{A,B} = \begin{cases}
(-q)^{d(A,B)} & \mathrm{if~}Y_{\mu^{-1}(\b)} \mathrm{~is~a~face~of~}A \\
0 & \mathrm{otherwise}.
\end{cases}\]

In the proof of the `numerical Koszulity condition'~\cite[Theorem 7.14]{KM}, this matrix $Y$ is proven to be invertible and so when we set $q=1$, it is still invertible.
\end{proof}

As we have shown that $\{\phi(A)\}$ forms a $\Z$-basis for $\cU$, the determinant of $S$ is equal to the determinant of the Brylawski--Schechtman--Varchenko form for the dual matroid (with the choice of weight function (or `collection of exponents') $a(p)=1$ for all $p \in I$).  By~\cite[Theorem 4.16]{BV} the determinant associated to the matroid $M$ is given by \[ \prod_{K \in \cF \setminus\{\emptyset\}} |K|^{\beta(K) \mu^+(M/K)}.\]
Applying the formula to the dual matroid gives \[ \det S(\cC) = \prod_{K \in \cF \setminus\{I\}} |I\setminus K|^{\beta(M/K) \mu^+(K)},\]
as we wished to show.

%!TEX root = ../main.tex

\section{Intersection matrices in hypertoric geometry and hypertoric category $\cO$}\label{sec-interpret}

In the special case where $\cC$ is a unimodular hyperplane arrangement, the non-$q$-version of the matrix $S$ has a geometric interpretation that appeared in the literature in~\cite{HS} where it is identified with a cohomological intersection form.  More precisely, if $\cC$ unimodular, then there is a corresponding smooth hypertoric variety $\tilde{\mathfrak{M}}$ defined by a hyperk\"ahler quotient construction.  The affinization of $\tilde{\mathfrak{M}}$ is a conical affine hypertoric variety $\mathfrak{M}$ and the fiber over the cone point is a union of toric varieties described by polytope regions of the arrangement $\cC$.  In this context, Hausel--Swartz show that $S$ arises as the intersection matrix on the top-dimensional cohomology of this fiber.

More generally, for any real hyperplane arrangement (and oriented matroid program) the matrix $S$ can also be viewed as arising from the Euler form on the Grothendieck group of hypertoric (or oriented matroid) category $\cO$.  Hypertoric category $\cO$ was introduced by Braden--Licata--Proudfoot--Webster~\cite{BLPWtorico} as a category of certain modules over a quantization deformation for hypertoric varieties, analogous to the Bernstein--Gelfand--Gelfand category $\cO$ of (possibly infinite-dimensional) representations of a semi-simple complex Lie algebra.  Moreover, they proved that their hypertoric category $\cO$ is equivalent to modules over a finite-dimensional quadratic algebra, defined from linear programming data and introduced in their paper~\cite{GDKD}.  In~\cite{KM}, Kowalenko and the second author generalized this notion from the setting of linear programming to that of oriented matroid programs. 

To be more precise, the input for defining hypertoric (resp. oriented matroid) category $\cO$ is a linear program (resp. oriented matroid program) meaning an essential real affine hyperplane arrangement $\cC$, generic as in the previous section (resp. an affine oriented matroid $(\tilde\cM,g)$, generic in an analogous way) together with a choice of generic covector $\xi$ (resp. generic single-element extension $(\tilde\cM,g,f)$ of the affine oriented matroid).  For any such linear program (resp. \textit{Euclidean} oriented matroid program), the resulting algebra $A=A(\cC,\xi)$ (resp. $A(\tilde\cM,f,g)$) is quasiheriditary and Koszul. 

The simple objects in the resulting category $\cO$, $A(\cC,\xi)-\mathrm{mod}$ (resp. $A(\tilde\cM,f,g)-\mathrm{mod}$) are parametrized by bases of the matroid underlying the hyperplane arrangement (resp. oriented matroid $\tilde{\cM}/g$) or correspondingly by regions of the arrangement that are bounded by $\xi$ (resp. topes that are bounded with respect to $f$).
For two compact regions/topes the number $S_q(A,B)$ arises in this context as the Euler form on the Grothendieck group between the corresponding simple modules.
\begin{theorem}\label{EulerForm} For a linear program (or Euclidean oriented matroid program) 
\[ S_q(A,B)=\langle [L_A], [L_B] \rangle_q=\sum_{i}(-q)^i\dim\operatorname{Ext}_{\cO}^i(L_A,L_B)\]
where $L_A,L_B\in \cO$ are the simple modules corresponding to the chambers $A,B$ bounded by $\xi$ (or $f$).
\end{theorem}
\begin{proof}
To prove this, we will pass to the Koszul dual side, where the corresponding statement is shown in \cite[Corollary 7.8]{KM}. Denote by $\cO^{\gr}=A-\mathrm{mod}^{\gr}$ the \emph{graded} category $\cO.$ We denote by $\langle n \rangle$ and $v:\cO^\Z\to \cO$ the functors shifting and forgetting the grading. 

We denote by $\cO^{\vee,\gr}$ the category associated to the Gale-dual linear/oriented matroid program. It arises as graded modules over the quadratic dual algebra $A^!.$ Then Koszul duality yields an equivalence of categories 
$$K:D^b(\cO^{\gr})\to D^b(\cO^{\vee,\gr}).$$
There is a natural bijection between chambers in the Gale dual programs and Koszul duality interchanges simple and projective modules as well as the graded and cohomological degrees
$$K(L_A\langle i\rangle)=P_A^{\vee}\langle -i\rangle[-i].$$
We hence obtain
\[
\begin{aligned}
    \sum_{i}(-q)^i&\dim\operatorname{Ext}_{\mathcal{O}}^i(L_A,L_B) \\
    &= \sum_{i}(-q)^i\dim\operatorname{Ext}_{\mathcal{O}^{\mathrm{gr}}}^i(L_A,L_B\langle i\rangle) \\
    &= \sum_{i}(-q)^i\dim\operatorname{Hom}_{\mathcal{O}^{\vee,\mathrm{gr}}}(P_A^\vee,P_B^\vee\langle -i\rangle) \\
   %&= \sum_{i}(-q)^i\dim (e_AA^!e_B)_i = \sum_{i}(-q)^i\dim (R_{A,B})_i \\
    &= (-q)^{d(A,B)} \sum_{i=0}^{d} \#\left\{\text{feasible vertices of } A\wedge B \text{ with } i \text{ outgoing edges}\right\} q^{2i} \\
    &= (-q)^{d(A,B)}h(A\wedge B)_{q^2}.
\end{aligned}
\]
The first equality uses the Koszulity of $\mathcal{O}^{\mathrm{gr}}$ which implies that $$\operatorname{Ext}_{\mathcal{O}^{\mathrm{gr}}}^i(L_A,L_B\langle j\rangle)=0$$ for all $j\neq i$. The second equality is a consequence of Koszul duality. The other equalities are shown in the proof of \cite[Corollary 7.8]{KM}.
\end{proof}

\appendix
%!TEX root = ../main.tex
\section{Examples of the $q$-determinant formula}\label{sec-examples}
\subsection{Points on a line}
Consider the arrangement $\mathcal{C}_n$ in $\R^1$, consisting of $n+1$ affine hyperplanes, say $H_i=\{-i\}$ for $i=1,\dots, n+1.$ We then obtain the bounded regions $A_i=(-i, -(i+1))$ for $i=1,\dots, n.$ So the intersection $\overline{A_i}\cap \overline{A_j}$ is either a line segment, a point or empty and we obtain
\begin{center}
\noindent
\begin{minipage}{0.5\textwidth}
\[
S(A_i,A_{j})=\begin{cases} 
    2 & \text{if } i = j, \\
    -1 & \text{if } |i-j|=1, \\
    0 & \text{else.}
\end{cases}
\]
\end{minipage}%
\begin{minipage}{0.5\textwidth}
\[
S_q(A_i,A_{j})=\begin{cases} 
    1+q^2 & \text{if } i = j, \\
    -q & \text{if } |i-j|=1, \\
    0 & \text{else.}
\end{cases}
\]
\end{minipage}
\end{center}
which yields the ($q$-)Schechtman--Varchenko matrices
\[
S(\mathcal{C}_n) = \begin{pmatrix}
    2 & -1 &  &  \\
    -1 & 2 & \ddots &  \\
     & \ddots & \ddots & -1 \\
     &  & -1 & 2
\end{pmatrix},
\quad
S_q(\mathcal{C}_n) = \begin{pmatrix}
    1+q^2 & -q &  &  \\
    -q & 1+q^2 & \ddots &  \\
     & \ddots & \ddots & -q \\
     &  & -q & 1+q^2
\end{pmatrix}.
\]
A straightforward induction shows that 
$$\det{S(\mathcal{C}_n)}=n+1 \text{ and }\det{S_q(\mathcal{C}_n)}=[n+1]_{q^2}.$$
In fact, $S(\mathcal{C}_n)$ is the Cartan matrix of the root system $A_n,$ and $S_q(\mathcal{C}_n)$ (up to a factor of $q$) the quantum Cartan matrix, see \cite{frenkelDeformationsWalgebrasAssociated1997}. 

On the other hand, the underlying matroid is isomorphic to the uniform matroid $U_{1,n+1}$ and $\emptyset$ is the only coloop-free set $\neq I.$ The product in Conjecture \ref{conj:mainconjecture} hence has a single factor
$$[|I\setminus\emptyset|]_{q^2}^{\beta(M)\mu^+(\emptyset)}=[n+1]_{q^2}^{\beta(U_{1,n+1})\mu^+(\emptyset)}=[n+1]_{q^2}.$$

\subsection{A Pseudoline Arrangement}
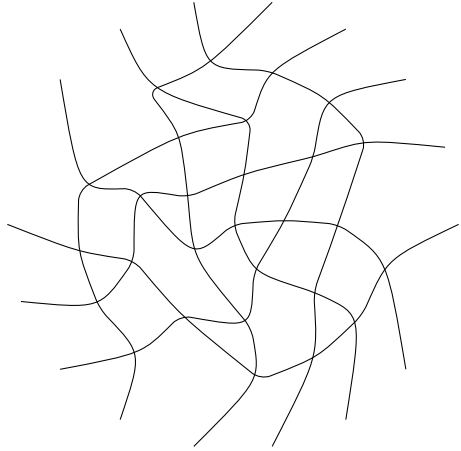
\begin{figure}[ht]
\begin{tikzpicture}
    \def\radius{3cm}
  \path (0:\radius) coordinate (2+);
  \path (20:\radius) coordinate (7+) ;
  \path (40:\radius) coordinate (8+) ;
  \path (60:\radius) coordinate  (6+) ;
  \path (80:\radius) coordinate (5+) ;
  \path (100:\radius) coordinate (9+) ;
  \path (120:\radius) coordinate (4+) ;
  \path (140:\radius) coordinate (3+) ;
 
  \path (180:\radius) coordinate (2-)  ;
  \path (200:\radius) coordinate (7-)  ;
  \path (220:\radius) coordinate (8-)  ;
  \path (240:\radius) coordinate (6-)   ;
  \path (250:\radius) coordinate (g);

  \path (260:\radius) coordinate (5-)   ;
  \path (280:\radius) coordinate (9-)   ;
  \path (300:\radius) coordinate (4-)   ;
  \path (320:\radius) coordinate (3-)   ;

  \coordinate (a) at (10:2.1cm);
  \coordinate (79) at (32.5:2.1cm);
  \coordinate (89) at (55:2.1cm);
  \coordinate (69) at (77.5:2.1cm);
  \coordinate (59) at (100:2.1cm);
  \coordinate (45) at (122.5:2.1cm);
  \coordinate (b) at (145:2.1cm);
  \coordinate (36) at (167.5:2.1cm);
  \coordinate (26) at (190:2.1cm);
  \coordinate (67) at (212.5:2.1cm);
  \coordinate (68) at (235:2.1cm);
  \coordinate (c) at (257.5:2.1cm);
  \coordinate (25) at (280:2.1cm);
  \coordinate (29) at (302.5:2.1cm);
  \coordinate (24) at (325:2.1cm);
  \coordinate (23) at (347.5:2.1cm);
  
  \coordinate (39) at (0:1.4cm);
  \coordinate (78) at (40:1.4cm);
  \coordinate (46) at (80:1.4cm);
  \coordinate (56) at (120:1.4cm);
  \coordinate (37) at (160:1.4cm);
  \coordinate (27) at (200:1.4cm);
  \coordinate (28) at (240:1.4cm);
  \coordinate (58) at (280:1.4cm);
  \coordinate (49) at (320:1.4cm);
  
  \coordinate (38) at (5:.7cm);
  \coordinate (47) at (77:.7cm);
  \coordinate (57) at (149:.7cm);
  \coordinate (35) at (221:.7cm);
  \coordinate (48) at (293:.7cm);
  
  \draw [rounded corners] (2+) .. controls (23) .. (24) .. controls (29) .. (25) .. controls (28) .. (27) .. controls (26) .. (2-);
  \draw [rounded corners] (7+) .. controls (79) .. (78) .. controls (47) .. (57) .. controls (37) .. (27) .. controls (67) .. (7-);
  \draw [rounded corners] (8+) .. controls (89) .. (78) .. controls (38) .. (48) .. controls (58) .. (28) .. controls (68) .. (8-);
  \draw [rounded corners] (6+) .. controls (69) .. (46) .. controls (56) .. (36) .. controls (26) .. (67) .. controls (68) .. (6-);
  \draw [rounded corners] (5+) .. controls (59) .. (45) .. controls (56) .. (57) .. controls (35) .. (58) .. controls (25) .. (5-);
  \draw [rounded corners] (9+) .. controls (59) .. (69) .. controls (89) .. (79) .. controls (39) .. (49) .. controls (29) .. (9-);
  \draw [rounded corners] (4+) .. controls (45) .. (46) .. controls (47) .. (0,0) .. controls (48) .. (49) .. controls (24) .. (4-); 
  \draw [rounded corners] (3+) .. controls (36) .. (37) .. controls (35) .. (0,0) .. controls (38) .. (39) .. controls (23) .. (3-); 
  
%   \draw[fill=gray, opacity=.2] (2-) arc [radius = \radius, start angle =  180,  end angle = 320] .. controls (23) .. (39) .. controls (79) .. (89) .. controls (69) .. (98:2.14cm) -- (45) -- (56) -- (36) -- (26) -- cycle;
  \end{tikzpicture}
  \caption{The pseudoline arrangement $\mathcal{C}_{\text{ps}}$}
  \label{fig:pseudoline}
\end{figure}
As a first example of Conjecture \ref{conj:om}, we consider the oriented matroid defined by the pseudoline arrangement $\mathcal{C}_{\text{ps}}$ depicted in Figure \ref{fig:pseudoline}, which has 21 compact topes. A computer calculation shows that 
$$\det(S_q(\mathcal{C}_{\text{ps}}))=[8]_{q^2}^6.$$
In this case the underlying (unoriented) matroid is the uniform matroid $U_{2,8}$.  The only coloop-free flats $K$ of $U_{2,8}$ are $\emptyset$ and $I$ itself.  Thus the product in the conjecture has a single factor, namely for $K=\emptyset$, and the conjectural formula reduces to 
\[ [| I\setminus \emptyset |]_{q^2}^{\beta(U_{2,8}/\emptyset) \mu^+(\emptyset)} =  [8]_{q^2}^{\beta(U_{2,8}) \mu^+(\emptyset)} = [8]_{q^2}^6,\]
where the last equality holds because $\mu^+(\emptyset) = 1$ and \[\beta(U_{2,8}) = \beta(U_{2,7})+\beta(U_{1,7}) =  \beta(U_{2,7}) + 1 = \beta(U_{2,6}) + 2 = \ldots = \beta(U_{2,2}) + 6 = 6.\]

\subsection{An affine V\'amos oriented matroid}

We conclude with an example where the underlying matroid is not representable.   Recall that the V\'amos matroid $V$ (rank 4 on 8 elements) is not representable over any field, but is orientable by~\cite{BLV}.  Let $\cV$ be the oriented matroid obtained as labelled and oriented in~\cite{BLV}. While there are many possible generic one-element lifts, for the purpose of checking an example, we simply pick one.  Namely, let $(\tilde{\cV},g)$ be the affine oriented matroid obtained from $\cV$  by adjoining $g$ to $\cV$ using the lexicographic one-element lift defined by the cobasis \{3, 6, 7, 8\}.  The set of feasible cocircuits for $(\tilde{\cV},g)$ is given in Table~\ref{cocircuits} and the corresponding set of bounded topes $\cP_0(\tilde\cV)$ is given in Table~\ref{topes}.  Using a computer calculation, we find that
\[ \det(S_q(\tilde\cV)) = [8]_{q^2}^{15} [4]_{q^2}^5.\]
On the other hand the V\'amos matroid $V$ has 5 proper nonempty coloop-free flats, namely the 5 circuits $K_1,\ldots,K_5$, each of size 4.  Computing the right hand side of Conjecture \ref{conj:mainconjecture} gives
\[\prod_{K \in \cF \setminus\{I\}} [| I\setminus K |]_{q^2}^{\beta(M/K) \mu^+(K)} = [| I |]_{q^2}^{\beta(V) \mu^+(\emptyset)} \cdot \prod_{i=1}^5 [| I\setminus K_i |]_{q^2}^{\beta(M/K_i) \mu^+(K_i)} \]
\[=[8]_{q^2}^{\beta(V) \mu^+(\emptyset)} \cdot ([4]_{q^2}^{\beta(U_{1,4}) \mu^+(U_{3,4})})^5 = [8]_{q^2}^{15} \cdot [4]_{q^2}^5.\]
Thus the conjecture holds for this example.

\begin{table}
\[\begin{array}{cccccccccc}
5 6 \overline{7} \overline{8} & 4 6 \overline{7} 8 &4 \overline{5} 7 8& 4 \overline{5} 6 8& 4 5 6 \overline{7}& \overline{3} \overline{6} \overline{7} \overline{8} & \overline{3} 5 \overline{7} \overline{8} &\overline{3} \overline{5} \overline{6} 8& \overline{3} \overline{5} \overline{6} \overline{7} &\overline{3} 4 \overline{7} 8 \\
 \overline{3} \overline{4} \overline{6} \overline{8} & \overline{3} 4 \overline{6} \overline{7}&\overline{3} 4 \overline{5} 8& \overline{3} 4 5 \overline{7}& \overline{3} \overline{4} \overline{5} \overline{6}& \overline{2} 6 \overline{7} 8& \overline{2} \overline{5} \overline{7} 8& \overline{2} \overline{5} \overline{6} 8& \overline{2} 5 6 \overline{7}& 2 4 6 8\\
 \overline{2} \overline{4} 6 \overline{7}& \overline{2} 4 \overline{5} 8& \overline{2} \overline{4} \overline{5} \overline{7}& \overline{2} \overline{3} \overline{7} 8& 2 \overline{3} \overline{6} \overline{8}& \overline{2} \overline{3} \overline{6} \overline{7}& 2 \overline{3} 5 \overline{8}& \overline{2} \overline{3} 5 \overline{7}& 2 \overline{3} \overline{5} \overline{6}& 2 \overline{3} 4 8\\
 \overline{2} \overline{3} \overline{4} \overline{7}& 2 \overline{3} 4 \overline{6}& 2 \overline{3} 4 5& 1 \overline{6} 7 8& 1 5 7 8&1 5 6 8& 1 5 6 \overline{7}& 1 4 7 8& 1 4 6 8& 1 \overline{4} \overline{6} 7\\
  \overline{1} 4 \overline{5} 8 &1 \overline{4} 5 7& 1 \overline{4} 5 6& 1 \overline{3} \overline{6} 8& 1 \overline{3} \overline{6} \overline{7}& 1 \overline{3} 5 8& 1 \overline{3} 5 \overline{7}&1 \overline{3} 4 8& \overline{1} \overline{3} 4 \overline{7}& 1 \overline{3} \overline{4} \overline{6} \\
  1 \overline{3} \overline{4} 5& 1 \overline{2} 7 8 & 1 \overline{2} 6 8&\overline{1} \overline{2} 6 \overline{7}& 1 \overline{2} \overline{5} 8& 1 2 5 7& 1 2 5 6& 1 \overline{2} \overline{4} 8& 1 \overline{2} \overline{4} \overline{7} & 1 \overline{2} \overline{4} \overline{6} \\
   1 \overline{2} \overline{4} 5& 1 \overline{2} \overline{3} 8& \overline{1} \overline{2} \overline{3} \overline{7}& 1 2 \overline{3} \overline{6}& 1 2 \overline{3} 5 & & & & &
  \end{array}\]
  \caption{Feasible cocircuits of $(\tilde{\cV},g)$}\label{cocircuits}
\end{table}

\begin{table}
\[
\begin{array}{ccccc}
1 2 \overline{3} \overline{4} \overline{5} \overline{6} \overline{7} \overline{8} &
1 \overline{2} \overline{3} 4 \overline{5} \overline{6} \overline{7} \overline{8}& 
1 2 \overline{3} 4 \overline{5} \overline{6} \overline{7} \overline{8}&
1 \overline{2} \overline{3} 4 5 \overline{6} \overline{7} \overline{8}&
1 2 \overline{3} 4 5 \overline{6} \overline{7} \overline{8} \\
1 \overline{2} \overline{3} 4 5 6 \overline{7} \overline{8}&
1 2 \overline{3} \overline{4} \overline{5} \overline{6} \overline{7} 8&
1 \overline{2} \overline{3} 4 \overline{5} \overline{6} \overline{7} 8&
1 2 \overline{3} 4 \overline{5} \overline{6} \overline{7} 8&
\overline{1} \overline{2} \overline{3} \overline{4} 5 \overline{6} \overline{7} 8 \\
1 \overline{2} \overline{3} \overline{4} 5 \overline{6} \overline{7} 8&
1 2 \overline{3} \overline{4} 5 \overline{6} \overline{7} 8&
\overline{1} \overline{2} \overline{3} 4 5 \overline{6} \overline{7} 8&
1 \overline{2} \overline{3} 4 5 \overline{6} \overline{7} 8&
1 2 \overline{3} 4 5 \overline{6} \overline{7} 8 \\
1 \overline{2} \overline{3} \overline{4} \overline{5} 6 \overline{7} 8&
1 \overline{2} \overline{3} 4 \overline{5} 6 \overline{7} 8&
1 2 \overline{3} 4 \overline{5} 6 \overline{7} 8&
\overline{1} \overline{2} \overline{3} \overline{4} 5 6 \overline{7} 8&
1 \overline{2} \overline{3} \overline{4} 5 6 \overline{7} 8 \\
1 2 \overline{3} \overline{4} 5 6 \overline{7} 8&
1 \overline{2} \overline{3} 4 5 6 \overline{7} 8&
1 \overline{2} \overline{3} \overline{4} 5 \overline{6} 7 8&
1 \overline{2} \overline{3} 4 5 \overline{6} 7 8&
1 \overline{2} \overline{3} \overline{4} \overline{5} 6 7 8 \\
\overline{1} \overline{2} \overline{3} 4 \overline{5} 6 7 8&
1 \overline{2} \overline{3} 4 \overline{5} 6 7 8&
1 \overline{2} \overline{3} \overline{4} 5 6 7 8&
1 2 \overline{3} \overline{4} 5 6 7 8&
1 \overline{2} \overline{3} 4 5 6 7 8
\end{array}\]
  \caption{Bounded topes of $(\tilde{\cV},g)$}\label{topes}
\end{table}

\bibliography{./refs.bib}
\bibliographystyle{amsalpha}

\end{document}